\RequirePackage{fix-cm}
\documentclass[twocolumn]{my-svjour3}      

\smartqed  

\usepackage{graphicx}
\usepackage{mathptmx} 
\usepackage{amssymb,amsmath}
\usepackage[misc,geometry]{ifsym} 
\usepackage{hyperref}

\newcommand{\doi}[1]{{\scriptsize 
 \textsc{doi}: \href{http://dx.doi.org/#1}{\nolinkurl{#1}}}}

\journalname{Soft Computing}

\begin{document}

\title{A new spectral method based on two classes of hat functions\\ 
for solving systems of fractional differential equations\\ 
and an application to respiratory syncytial virus infection}

\titlerunning{A new spectral method for solving systems of FDEs and an application to RSV}

\author{Somayeh Nemati \and Delfim F. M. Torres}

\authorrunning{S. Nemati \and D. F. M. Torres} 

\institute{S. Nemati \at
Department of Mathematics,\\ 
Faculty of Mathematical Sciences,\\ 
University of Mazandaran, Babolsar, Iran\\
\email{s.nemati@umz.ac.ir}\\
\emph{Present address:} 
Center for Research and Development in Mathematics and Applications (CIDMA),
Department of Mathematics, University of Aveiro, 3810-193 Aveiro, Portugal\\
\email{s.nemati@ua.pt}\\
ORCID: \url{https://orcid.org/0000-0003-1724-6296}
\and
D. F. M. Torres \Letter \at
Center for Research and Development in Mathematics and Applications (CIDMA), 
Department of Mathematics, University of Aveiro, 3810-193 Aveiro, Portugal\\
Tel.: +351-234-370668\\
Fax: +351-234-370066\\
\email{delfim@ua.pt} \\
ORCID: \url{https://orcid.org/0000-0001-8641-2505}}


\date{Submitted: 26-April-2019; Revised: 02-Nov-2019; Accepted: 17-Dec-2019}

\maketitle


\begin{abstract}
We propose a new spectral method, based on two classes of hat functions, 
for solving systems of fractional differential equations. The fractional 
derivative is considered in the Caputo sense. Properties of the basis functions, 
Caputo derivatives, and Riemann--Liouville fractional integrals, are used 
to reduce the main problem to a system of nonlinear algebraic equations. 
By analyzing in detail the resulting system, we show that the method needs 
few computational efforts. Two test problems are considered to illustrate 
the efficiency and accuracy of the proposed method. Finally, an application 
to a recent mathematical model in epidemiology is given, precisely 
to a system of fractional differential equations modeling
the respiratory syncytial virus infection.
\end{abstract}


\section{Introduction}
\label{sec:1}

In the last few decades, it has been shown that fractional calculus appears 
in modeling of many real-world phenomena, in different branches of science 
where nonlocality plays an essential role, such as physics, chemistry, biology, 
economics, engineering, signal and image processing, and control theory 
\cite{Agarwal,Baleanu},\\
\cite{Malinowska,Morales},\\
\cite{Rekhviashvili,Tariboon}. 
This is mainly due to the fact that fractional 
operators consider the evolution of a system, by taking the global correlation, 
and not only local characteristics \cite{Almeida}. As a result, fractional 
systems have attracted the attention of many researchers in different areas. 
Nevertheless, obtaining analytic solutions for such problems is very difficult. 
So, in most cases, the exact solution is not known, and it is necessary 
to find a numerical approximation. Therefore, many researchers have worked 
on numerical methods to obtain some numerical solutions of fractional 
dynamic systems (see, e.g., \cite{MR3904404,Jafarian1,Jafarian2,El-Sayed,Nigmatullin}).
 
Hat basis functions consist of a set of piecewise continuous functions with shape of hats, 
when plotted in two dimensional planes. These functions are usually defined on the 
interval $[0,1]$ and a generalization of them is obtained by extending the interval 
into $[0,\tau]$ with any arbitrary positive number $\tau$. Hat basis functions have 
shown to be a powerful mathematical tool in solving many different kinds of equations. 
For example, in \cite{Babolian}, hat functions have been used to solve linear and 
nonlinear integral equations of second kind. These functions have been also
efficiently employed to solve fractional differential equations (FDEs) \cite{Tripathi},
while in \cite{Heydari} generalized hat functions are used for solving stochastic 
It\^{o}--Volterra integral equations. Recently, a modification 
of hat functions has been introduced and used in order to solve a variety of problems. 
To mention some of these problems, we refer to two-dimensional linear Fredholm integral 
equations \cite{Mirzaee1}, integral equations of Stratonovich--Volterra \cite{Mirzaee2} 
and Volterra--Fredholm type \cite{Mirzaee3}, and fractional integro-differential 
\cite{Nemati1}, fractional pantograph nonlinear differential equations\\ 
\cite{Nemati2} and fractional optimal control problems \cite{Nemati3}. 
The aim of our work is to have a comparison between these two classes 
of hat basis functions in solving systems of FDEs. 
 
The paper is organized as follows. Section~\ref{sec:2} is devoted to the 
required preliminaries for presenting our numerical technique. In Section~\ref{sec:3}, 
we present the new numerical method, which is based on two classes of hat functions 
for solving systems of FDEs. Section~\ref{sec:4} is concerned with an application 
of the method for solving a problem in epidemiology related to  human respiratory 
syncytial virus. Finally, concluding remarks are given in Section~\ref{sec:5}.  


\section{Preliminaries}
\label{sec:2}

In this section, some necessary definitions and properties of fractional calculus 
are presented. Moreover, two classes of hat functions and some 
of their properties are recalled.


\subsection{Preliminaries of fractional calculus}

In this work, we employ fractional differentiation in the sense 
of Caputo, which is defined via the Riemann--Liouville fractional integral.

\begin{definition} \cite{Podlubny}
The (left) Riemann--Liouville fractional integral operator with order 
$\alpha\geq 0$ of a given function $y$ is defined as 
\begin{equation*}
\begin{split}
&_{0}I_t^{\alpha}y(t)=\frac{1}{\Gamma(\alpha)}\int_0^t(t-s)^{\alpha-1}y(s)ds,\\
&_{0}I_t^0y(t)=y(t),
\end{split}
\end{equation*}
where $\Gamma(\cdot)$ is the Euler gamma function.
\end{definition}

\begin{definition} \cite{Podlubny}
The (left) Caputo fractional derivative of order $\alpha>0$ of a function $y$ is defined as
\begin{equation*}
^C_0D_t^{\alpha}y(t)=\frac{1}{\Gamma(m-\alpha)}
\int_0^t(t-s)^{m-\alpha-1}y^{(m)}(s)ds,
\end{equation*}
where $m-1<\alpha \leq m$.
\end{definition}

For $m-1<\alpha\leq m$, $m\in\mathbb{N}$, we recall two important properties 
of the Caputo derivative and Riemann--Liouville integral:
\begin{equation*}
^C_0D_t^\alpha (_{0}I^\alpha y(t))=y(t),
\end{equation*}
\begin{equation}
\label{2.1}
_{0}I_t^\alpha ({^C_0D_t^\alpha} y(t))
=y(t)-\sum_{i=0}^{m-1}y^{(i)}(0)\frac{t^{i}}{i!},
\quad t>0.
\end{equation}


\subsection{Hat functions}

We consider both generalized and modified hat functions.

\noindent {\it Generalized hat functions.}
The interval $[0,\tau]$ is divided into $n$ subintervals $[ih,(i+1)h]$, 
$i = 0,1,2,\ldots,n-1$, of equal length $h$, where $h=\frac{\tau}{n}$. 
Then, the generalized hat functions (GHFs) are defined as follows \cite{Tripathi}:
\begin{equation*}
\psi^{G}_0(t)=\left\{
\begin{array}{ll}
\frac{h-t}{h},&0\leq t\leq h,\\
&\\
0,&\text{otherwise},
\end{array}\right.
\end{equation*}
\begin{equation*}
\psi^G_i(t)
=\left\{
\begin{array}{ll}
\frac{t-(i-1)h}{h},&(i-1)h\leq t\leq ih,\\
&\\
\frac{(i+1)h-t}{h},&ih\leq t\leq(i+1)h,\\
&\\
0,&\text{otherwise},
\end{array}\right.
\end{equation*}
\begin{equation*}
\psi^G_n(t)=\left\{
\begin{array}{ll}
\frac{t-(\tau-h)}{h},&\tau-h\leq t\leq \tau,\\
&\\
0,&\text{otherwise}.
\end{array}\right.
\end{equation*}
These functions form a set of $(n+1)$ linearly independent continuous 
functions in $L^2[0,\tau]$, satisfying the property
\begin{equation}
\label{2.1.1}
\psi_i^G(jh)=\left\{\begin{array}{ll}1,&i=j,\\
&\\
0,&i\neq j.
\end{array}\right.
\end{equation}
An arbitrary function $y\in L^2[0,\tau]$ can be approximated using the GHFs as follows:
\begin{equation}
\label{2.2}
y(t)\simeq y_n(t)=\sum_{i=0}^n y(ih) \psi^G_i(t)=A^T \Psi_G(t),
\end{equation}
where
\begin{equation}
\label{2.3}
\Psi_G(t)=[\psi^G_0(t),\psi^G_1(t),\ldots,\psi^G_n(t)]^T
\end{equation}
and
\begin{equation*}
A=[a_0,a_1,\ldots,a_n]^T
\end{equation*}
with $a_i=y(ih)$.

\begin{theorem} \cite{Tripathi}
If $y\in C^2([0,\tau])$ is approximated 
by the family of first $(n+1)$ GHFs as \eqref{2.2}, then
\begin{equation*}
|y(t)-y_n(t)|=O(h^2).
\end{equation*}
\end{theorem}

Let $\Psi_G$ be the GHFs basis vector given by \eqref{2.3} and $\alpha>0$. Then,
\begin{equation}
\label{2.4}
I_t^\alpha\Psi_G(t)\simeq P_G^{\alpha}\Psi_G(t),
\end{equation}
where $P_G^{\alpha}$ is a matrix of dimension $(n+1)\times(n+1)$  
called the operational matrix of fractional integration of order $\alpha$ of the GHFs. 
This matrix is given as \cite{Tripathi}  
\begin{equation}
\label{2.5}
P_G^{\alpha}
=\frac{h^\alpha}{\Gamma(\alpha+2)}\left[
\begin{array}{cccccc}
0 & \zeta_1  &\zeta_2& \zeta_3  &\ldots &\zeta_n \\
0 &1 & \varrho_1 & \varrho_2&\ldots & \varrho_{n-1} \\
0 & 0 &1 &\varrho_1&\ldots & \varrho_{n-2}\\
0 & 0 & 0& 1 &  \ldots &\varrho_{n-3}\\
\ldots & \ldots &\ldots& \ldots  & \ldots & \ldots\\
0 & 0 & 0 &0 & \ldots & 1\\
\end{array}
\right],
\end{equation}
where
\begin{equation*}
\zeta_i=i^\alpha(\alpha-i+1)+(i-1)^{\alpha+1}, 
\quad i=1,2,\ldots,n,
\end{equation*}
and
\begin{equation*}
\varrho_i=(i+1)^{\alpha+1}-2i^{\alpha+1}+(i-1)^{\alpha+1}, \quad i=1,2,\ldots,n-1.\\
\end{equation*}


\noindent {\it Modified hat functions.}
By considering an even integer number $n$, the interval $[0,\tau]$ is divided 
into $n$ subintervals $[ih,(i+1)h]$, $i=0,1,2,\ldots,n-1$, with equal length 
$h=\frac{\tau}{n}$. The modified hat functions (MHFs) form a set of $(n+1)$ 
linearly independent functions in $L^2[0,\tau]$. These functions are defined 
as follows \cite{Nemati1,Nemati2}:
\begin{equation*}
\psi^{M}_0(t)=\left\{
\begin{array}{ll}
\frac{1}{2h^2}(t-h)(t-2h),&0\leq t\leq 2h,\\
&\\
0,&\text{otherwise};
\end{array}\right.
\end{equation*}
if $i$ is odd and $1\leq i\leq n-1$, then
\begin{equation*}
\psi^M_i(t)=\left\{
\begin{array}{ll}
\frac{-(t-(i-1)h)}{h^2}(t-(i+1)h),&(i-1)h\leq t\leq(i+1)h,\\
&\\
0,&\text{otherwise};
\end{array}\right.
\end{equation*}
if $i$ is even and $2\leq i\leq n-2$, then
\begin{equation*}
\psi^M_i(t)
=\left\{
\begin{array}{ll}
\frac{1}{2h^2}(t-(i-1)h)(t-(i-2)h),&(i-2)h\leq t\leq ih,\\
&\\
\frac{1}{2h^2}(t-(i+1)h)(t-(i+2)h),&ih\leq t\leq(i+2)h,\\
&\\
0,&\text{otherwise},
\end{array}\right.
\end{equation*}
and
\begin{equation*}
\psi^M_n(t)=\left\{
\begin{array}{ll}
\frac{1}{2h^2}(t-(\tau-h))(t-(\tau-2h)),&\tau-2h\leq t\leq \tau,\\
&\\
0,&\text{otherwise}.
\end{array}\right.
\end{equation*}
The following property is satisfied for the set of MHFs:
\begin{equation}
\label{2.1.2}
\psi_i^M(jh)=\left\{\begin{array}{ll}1,
&i=j,\\
&\\
0,&i\neq j.
\end{array}\right.
\end{equation}
Any function $y\in L^2[0,\tau]$ may be approximated using the MHFs as
\begin{equation}
\label{2.6}
y(t)\simeq y_n(t)=\sum_{i=0}^n y(ih) \psi^M_i(t)=A^T \Psi_M(t),
\end{equation}
where
\begin{equation}
\label{2.7}
\Psi_M(t)=[\psi^M_0(t),\psi^M_1(t),\ldots,\psi^M_n(t)]^T
\end{equation}
and
\begin{equation*}
A=[a_0,a_1,\ldots,a_n]^T
\end{equation*}
with $a_i=y(ih)$.

\begin{theorem} \cite{Nemati1}
If $y\in C^3([0,\tau])$ is approximated by the family 
of first $(n + 1)$ MHFs as \eqref{2.6}, then
\begin{equation*}
|y(t)-y_n(t)|=O(h^3).
\end{equation*}
\end{theorem}

Let $\Psi_M$ be the MHFs basis vector given by \eqref{2.7} and $\alpha>0$. Then,
\begin{equation}
\label{2.8}
I_t^\alpha\Psi_M(t)\simeq P_M^{\alpha}\Psi_M(t),
\end{equation}
where $P_M^{\alpha}$ is a matrix of dimension $(n+1)\times(n+1)$ called 
the operational matrix of fractional integration of order $\alpha$ 
of the MHFs. This matrix is given (see \cite{Nemati1,Nemati2}) as  
\begin{equation}
\label{2.9}
P_M^{\alpha}
=\frac{h^\alpha}{2\Gamma(\alpha+3)}\left[
\begin{array}{cccccccc}
0 & \beta_1  &\beta_2& \beta_3 &\beta_4 &\ldots & \beta_{n-1}& \beta_n \\
0 & \eta_0 & \eta_1 & \eta_2&\eta_3&\ldots &\eta_{n-2} & \eta_{n-1} \\
0 & \xi_{-1} &\xi_0 & \xi_1&\xi_2&\ldots &\xi_{n-3}& \xi_{n-2}\\
0 & 0 & 0 & \eta_0 & \eta_1 & \ldots &\eta_{n-4} & \eta_{n-3}\\
0 & 0 & 0 &\xi_{-1} &\xi_0 & \ldots &\xi_{n-5} & \xi_{n-4}\\
\ldots & \ldots &\ldots& \ldots& \ldots  & \ldots &\ldots  & \ldots\\
0 & 0 & 0 & 0&0&\ldots &\eta_0 & \eta_1 \\
0 & 0 & 0 &0 &0 & \ldots &\xi_{-1} & \xi_0\\
\end{array}
\right]
\end{equation}
with
\begin{equation*}
\begin{split}
&\beta_1=\alpha(3+2\alpha),\\
&\beta_i=i^{\alpha+1}(2i-6-3\alpha)+2i^\alpha(1+\alpha)(2+\alpha)\\
&\qquad -(i-2)^{\alpha+1}(2i-2+\alpha), \quad i=2,3,\ldots,n,\\
&\eta_0=4(1+\alpha),\\
&\eta_i=4[(i-1)^{\alpha+1}(i+1+\alpha)-(i+1)^{\alpha+1}(i-1-\alpha)],\\
&\quad i=1,2,\ldots,n-1,\\
&\xi_{-1}=-\alpha,\\
&\xi_0=2^{\alpha+1}(2-\alpha),\\
&\xi_1=3^{\alpha+1}(4-\alpha)-6(2+\alpha),\\
&\xi_i=(i+2)^{\alpha+1}(2i+2-\alpha)-6i^{\alpha+1}(2+\alpha)\\
&\qquad -(i-2)^{\alpha+1}(2i-2+\alpha), \quad i=2,3,\ldots,n-2.
\end{split}
\end{equation*}


\section{Main results}
\label{sec:3}

We consider the following general initial value problem, 
described by a system of $m$ FDEs
of order $0<\alpha\leq 1$:
\begin{equation}
\label{3.1}
\begin{cases}
{_0^C D_t^\alpha} y_1(t)=f_1(t,y_1(t),\ldots,y_m(t)),& y_1(0)=y_{1,0}\\
~~~\vdots\\
{_0^C D_t^\alpha} y_m(t)=f_m(t,y_1(t),\ldots,y_m(t)),& y_m(0)=y_{m,0}.
\end{cases}
\end{equation}
The aim is to seek functions $y_1$, \ldots, $y_m$ 
solution of \eqref{3.1} on the interval $0\leq t\leq\tau$. 


\subsection{Numerical method}

In order to find a numerical solution of \eqref{3.1}, we consider approximations 
of the fractional derivative of the unknown functions as follows:
\begin{equation}
\label{3.3}
{_0^C D_t^\alpha} y_1(t)=A_1^T\Psi(t),\ldots,{_0^C D_t^\alpha} y_m(t)=A_m^T\Psi(t),
\end{equation}
where $A_i$ are the coefficients vectors with the unknown elements $a_{i,j}$, 
$j=0,1,\ldots,n$, 
\begin{equation}
A_i=[a_{i,0},a_{i,1},\ldots,a_{i,n}]^T,\quad i=1,\ldots,m,
\end{equation}
and $\Psi(t)$ equals to either $\Psi_G(t)$, corresponding to the GHFs 
as the basis functions, or $\Psi_M(t)$, corresponding to the MHFs 
basis functions. Then, using \eqref{2.1} and the initial conditions 
given in \eqref{3.1}, we can write
\begin{equation*}
y_i(t)= A_i^T I_t^\alpha \Psi(t)+y_{i,0}.
\end{equation*}
Using \eqref{2.4} or \eqref{2.8}, according to the chosen basis functions, we get
\begin{equation}
\label{3.2}
y_i(t)=A_i^TP^{\alpha}\Psi(t)+y_{i,0}
=(A_i^TP^{\alpha}+d_i^T)\Psi(t)=Y_i^T\Psi(t),\\
\end{equation}
where
\begin{equation}
\label{3.6}
d_i=[y_{i,0},y_{i,0},\ldots,y_{i,0}]^T, 
\quad Y_i=(A_i^TP^{\alpha}+d_i^T)^T
=[Y_{i,0},\ldots,Y_{i,n}]^T
\end{equation}
and $P^{\alpha}$ is given by \eqref{2.5} for the GHFs, or \eqref{2.9} for the MHFs. 
Therefore, by setting $t=jh$ in \eqref{3.2} and employing \eqref{2.1.1} or \eqref{2.1.2}, 
we obtain $y_i(jh)=Y_{i,j}$. Now, we can obtain approximations of the functions 
$f_i:\mathbb{R}^{m+1}\rightarrow \mathbb{R}$, by using the considered basis functions, 
as follows:
\begin{multline}
\label{3.4}
f_i(t,y_1(t),\ldots,y_m(t))
=\sum_{j=0}^nf_i(jh,y_1(jh),\ldots,y_m(jh))\psi_j(t)\\
=\sum_{j=0}^nf_i(jh,Y_{1,j},\ldots,Y_{m,j})\psi_j(t)
=F_i(\Theta,Y_1,\ldots,Y_m)\Psi(t),
\end{multline}
where $\psi_j(t)=\psi_j^G(t)$ or $\psi_j^M(t)$, and
\begin{multline*}
F_i(\Theta,Y_1,\ldots,Y_m)=[f_i(0,Y_{1,0},\ldots,Y_{m,0}),\\
f_i(h,Y_{1,1},\ldots,Y_{m,1}),\ldots,f_i(\tau,Y_{1,n},\ldots,Y_{m,n})].
\end{multline*}
By substitution \eqref{3.3} and \eqref{3.4} into \eqref{3.1}, we have
\begin{equation*}
\begin{split}
&A_1^T\Psi(t)=F_1(\Theta,Y_1,\ldots,Y_m)\Psi(t),\\
&~~~\vdots\\
&A_m^T\Psi(t)=F_m(\Theta,Y_1,\ldots,Y_m)\Psi(t),\\
\end{split}
\end{equation*}
which gives the following system:
\begin{equation*}
\begin{split}
&A_1^T-F_1(\Theta,Y_1,\ldots,Y_m)=0,\\
&~~~\vdots\\
&A_m^T-F_m(\Theta,Y_1,\ldots,Y_m)=0,\\
\end{split}
\end{equation*}
or
\begin{equation}\label{3.5}
a_{i,j}=f_i(jh,Y_{1,j},\ldots,Y_{m,j}),\quad i=1,\ldots,m,~j=0,\ldots,n.
\end{equation}
This system includes $m(n+1)$ nonlinear algebraic equations with $m(n+1)$ 
unknown parameters, which are the elements of $A_i$, $i=1,\ldots,m$. 
By solving this system, approximations of the functions 
$y_i$ are given by \eqref{3.2}.


\subsection{Complexity of the resulting systems}

The speed of the numerical method given above depends on the speed 
of solving system \eqref{3.5}. Therefore, the form of this system 
is an important aspect of our method. Here, we display the form 
of the system given in \eqref{3.5} in detail, 
for each of the basis functions.

\noindent {\bf Case 1:} {\it GHFs}.

\noindent According to the form of the operational matrix of fractional 
integration of the GHFs, given by \eqref{2.5}, 
we can rewrite this matrix as follows:
\begin{equation*}
P_G^{\alpha}=\left[p_{i,j}^G\right],\quad i,j=0,1,\ldots,n,
\end{equation*}
with
\begin{equation*}
\begin{split}
&p_{i,0}^G=0,\quad \text{for $i=0,1,2,\ldots,n$},\\
&p_{i,j}^G=0,\quad \text{for $1\leq j<i\leq n$}.\\
\end{split}
\end{equation*}
Therefore, the elements of the vectors $Y_i$, 
$i=1,\ldots,m$, in \eqref{3.2} are given, 
using \eqref{3.6}, by
\begin{equation}
\label{3.7}
\begin{split}
&Y_{i,0}=y_{i,0},\\
&Y_{i,j}=\sum_{k=0}^ja_{i,k}p_{k,j}^G+y_{i,0},\quad j=1,\ldots,n.
\end{split}
\end{equation}
Taking \eqref{3.7} into account, we rewrite system \eqref{3.5} 
as follows:
\begin{subequations}
\begin{align}
&a_{i,0}=f_i(0,y_{1,0},\ldots,y_{m,0}),\label{3.8.1}\\
&a_{i,1}=f_i(h,\sum_{k=0}^1a_{1,k}p_{k,1}^G+y_{1,0},
\ldots,\sum_{k=0}^1a_{m,k}p_{k,1}^G+y_{m,0}),\label{3.8.2}\\
&~~~\vdots\nonumber\\
&a_{i,n}=f_i(\tau,\sum_{k=0}^na_{1,k}p_{k,n}^G+y_{1,0},
\ldots,\sum_{k=0}^na_{m,k}p_{k,n}^G+y_{m,0}),\label{3.8.3}
\end{align}
\end{subequations}
$i=1,\ldots,m$. As it can be seen in \eqref{3.8.1}, the values 
of the unknown parameters $a_{i,0}$, $i=1,\ldots,m$, are obtained 
easily by using the initial conditions. By substituting the given 
$a_{i,0}$ into \eqref{3.8.2}, we have a system of $m$ nonlinear 
algebraic equations in unknown parameters $a_{i,1}$, $i=1,\ldots,m$. 
After solving this system, and substituting the obtained results 
for $a_{i,1}$ into the $m$ next equations, a system of $m$ equations 
in $a_{i,2}$, $i=1,\ldots,m$, is obtained. This process continues 
until $a_{i,n}$, $i=1,\ldots,m$, are found by solving \eqref{3.8.3}, 
in which the results for $a_{i,0}$, $a_{i,1}$, \ldots, $a_{i,n-1}$ 
have been substituted. Therefore, our method, based on GHFs, reduces 
the main problem to the solution of $n$ systems of $m$ 
nonlinear algebraic equations.

\noindent {\bf Case 2:} {\it MHFs}.

\noindent In a similar way as for GHFs, we rewrite the operational 
matrix of fractional integration of MHFs as follows:
\begin{equation*}
P_M^{\alpha}=\left[p_{i,j}^M\right],
\quad i,j=0,1,\ldots,n,
\end{equation*}
with
\begin{equation}
\label{3.9}
\begin{split}
&p_{i,0}^M=0,\quad \text{for $i=0,1,2,\ldots,n$},\\
&p_{i,j}^M=0,\quad \text{for $j=1,3,\ldots,n-1,$ $i=j+2,\ldots,n$},\\
&p_{i,j}^M=0,\quad \text{for $j=2,4,\ldots,n,$ $i=j+1,\ldots,n$}.
\end{split}
\end{equation}
By considering \eqref{3.9} for writing the elements of the vectors 
$Y_i$, $i=1,\ldots,m$, in \eqref{3.2}, we get
\begin{equation}
\label{3.10}
\begin{split}
&Y_{i,0}=y_{i,0},\\
&Y_{i,j}=\sum_{k=0}^{j+1}a_{i,k}p_{k,j}^M+y_{i,0},\quad j=1,3,\ldots,n-1\\
&Y_{i,j}=\sum_{k=0}^{j}a_{i,k}p_{k,j}^M+y_{i,0},\quad j=2,4,\ldots,n.
\end{split}
\end{equation}
By substituting the values of $Y_{i,j}$ given in \eqref{3.10} into system 
\eqref{3.5}, this system can be rewritten in detail as follows:
\begin{subequations}
\begin{align}
&a_{i,0}=f_i\left(0,y_{1,0},\ldots,y_{m,0}\right),\label{3.11.1}\\
&a_{i,1}=f_i\left(h,\sum_{k=0}^2a_{1,k}p_{k,1}^M+y_{1,0},\ldots,
\sum_{k=0}^2a_{m,k}p_{k,1}^M+y_{m,0}\right),\label{3.11.2}\\
&a_{i,2}=f_i\left(2h,\sum_{k=0}^2a_{1,k}p_{k,2}^M+y_{1,0},\ldots,
\sum_{k=0}^2a_{m,k}p_{k,2}^M+y_{m,0}\right),\label{3.11.3}\\
&~~~\vdots\nonumber\\
&a_{i,n-1}=f_i\left(\tau-h,\sum_{k=0}^na_{1,k}p_{k,n-1}^M+y_{1,0},
\ldots,\sum_{k=0}^na_{m,k}p_{k,n-1}^M+y_{m,0}\right),\label{3.11.4}\\
&a_{i,n}=f_i\left(\tau,\sum_{k=0}^na_{1,k}p_{k,n}^M+y_{1,0},
\ldots,\sum_{k=0}^na_{m,k}p_{k,n}^M+y_{m,0}\right),\label{3.11.5}
\end{align}
\end{subequations}
$i=1,\ldots,m$. It is seen that the values of the unknown parameters 
$a_{i,0}$, $i=1,\ldots,m$, are given easily by substituting the initial 
conditions into \eqref{3.11.1}. By substituting the given results of 
$a_{i,0}$ into \eqref{3.11.2} and \eqref{3.11.3}, a system of $2m$ 
nonlinear algebraic equations in the unknown parameters $a_{i,1}$ 
and $a_{i,2}$, $i=1,\ldots,m$, is obtained. After solving this system, 
and substituting the obtained results for $a_{i,1}$ and $a_{i,2}$ 
into the $2m$ next equations, a system of $2m$ algebraic equations 
in $a_{i,3}$ and $a_{i,4}$, $i=1,\ldots,m$, is given. By continuation 
of this process, we find $a_{i,n-1}$ and $a_{i,n}$, $i=1,\ldots,m$, 
by solving \eqref{3.11.4} and \eqref{3.11.5}, in which the results 
for $a_{i,0}$, $a_{i,1}$, \ldots, $a_{i,n-2}$ have been substituted. 
Our method based on MHFs reduces the system of $m$ 
nonlinear FDEs to solving $\frac{n}{2}$ systems of $2m$ 
nonlinear algebraic equations.


\subsection{Test problems}

In order to illustrate the efficiency and accuracy of the proposed method, 
we apply it to two test problems whose exact solutions are known. 
\begin{example}\label{ex1}
Consider the following system of FDEs:
\begin{equation}
\label{3.12}
\begin{cases}
{_0^CD_t^{0.5}}y_1(t)=\sqrt{t}y_1(t)-y_2(t)
+\frac{15\sqrt{\pi}}{16}t^2,\quad   y_1(0)=0,\\
{_0^CD_t^{0.5}}y_2(t)=\frac{16}{5\sqrt{\pi}}y_1(t)
+y_2^2(t)-t^6,\quad \quad \quad y_2(0)=0,
\end{cases}
\end{equation}
on the interval $0\leq t \leq 1$, which has the exact solution
\begin{equation}
y_1(t)=t^{2.5},\quad y_2(t)=t^3.
\end{equation}
We have solved this problem by the suggested method, based 
on GHFs and MHFs, with different values of $n$, and display 
the numerical results in Figure~\ref{fig:1}, Table~\ref{tab:1}, 
and Figure~\ref{fig:2}. In Figure~\ref{fig:1}, the approximate 
solutions of $y_1(t)$ and $y_2(t)$, obtained by employing our 
method with $n=2$, are shown. In this figure, the numerical results 
of $y_1(t)$ and $y_2(t)$, obtained by GHFs, are, respectively, displayed 
by $y_1^G(t)$ and $y_2^G(t)$. Also, the numerical results of the unknown 
functions given by MHFs are represented by $y_1^M(t)$ and $y_2^M(t)$. 
In Table~\ref{tab:1}, we see the numerical results for the functions 
$y_1(t)$ and $y_2(t)$, with different values of $n$, together with the 
CPU time (in seconds), which have been obtained on a 2.5 GHz Core i7 
personal computer with 16~GB of RAM using \textsf{Mathematica 11.3}. 
For solving the resulting systems of algebraic equations, the \textsf{Mathematica}
function \textsf{FindRoot} was used. In this table, the following notations 
are used for introducing the error and the convergence order of the method:
\begin{equation*}
\begin{split}
&e_{1,n}^G=\max_{0\leq i\leq n}|y_1(ih)-y_{1,n}^G(ih)|,
\quad \rho_{1,n}^G=\log_2{\left(\frac{e_{1,n}^G}{e_{1,2n}^G}\right)},\\
&e_{1,n}^M=\max_{0\leq i\leq n}|y_1(ih)-y_{1,n}^M(ih)|,
\quad \rho_{1,n}^M=\log_2{\left(\frac{e_{1,n}^M}{e_{1,2n}^M}\right)},\\
&e_{2,n}^G=\max_{0\leq i\leq n}|y_2(ih)-y_{2,n}^G(ih)|,
\quad \rho_{2,n}^G=\log_2{\left(\frac{e_{2,n}^G}{e_{2,2n}^G}\right)},\\
&e_{2,n}^M=\max_{0\leq i\leq n}|y_2(ih)-y_{2,n}^M(ih)|,
\quad \rho_{2,n}^M=\log_2{\left(\frac{e_{2,n}^M}{e_{2,2n}^M}\right)},
\end{split}
\end{equation*}
where $y_1(t)$ and $y_2(t)$ are the exact solutions, $y_{1,n}^G(t)$ 
and $y_{2,n}^G(t)$ are the approximate solutions obtained with GHFs, 
and $y_{1,n}^M(t)$ and $y_{2,n}^M(t)$ are the approximate solutions 
obtained with MHFs. These results confirm the $O(h^2)$ accuracy order 
of the numerical method with GHFs and the $O(h^3)$ accuracy order 
of the numerical method with MHFs. Finally, in Figure~\ref{fig:2} (left), 
the results for the errors obtained by employing our method, for some selected 
values of $n$, are plotted in a logarithmic scale. Moreover, the CPU times 
of the method are plotted in Figure \ref{fig:2} (right). It can be seen that 
the computational complexity of the resulting systems of GHFs and MHFs are similar. 
\begin{figure*}
\centering
\includegraphics[scale=0.70]{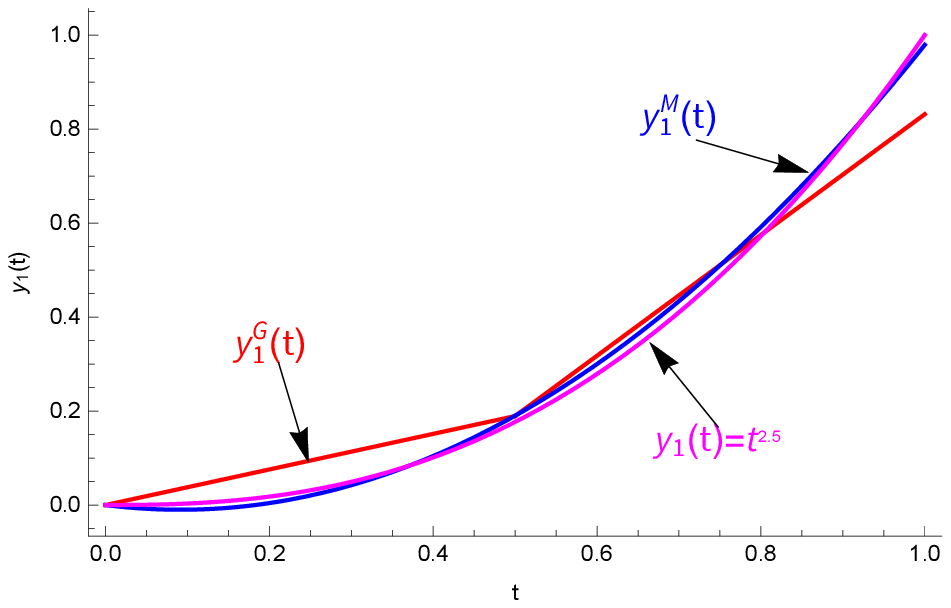}
\includegraphics[scale=0.70]{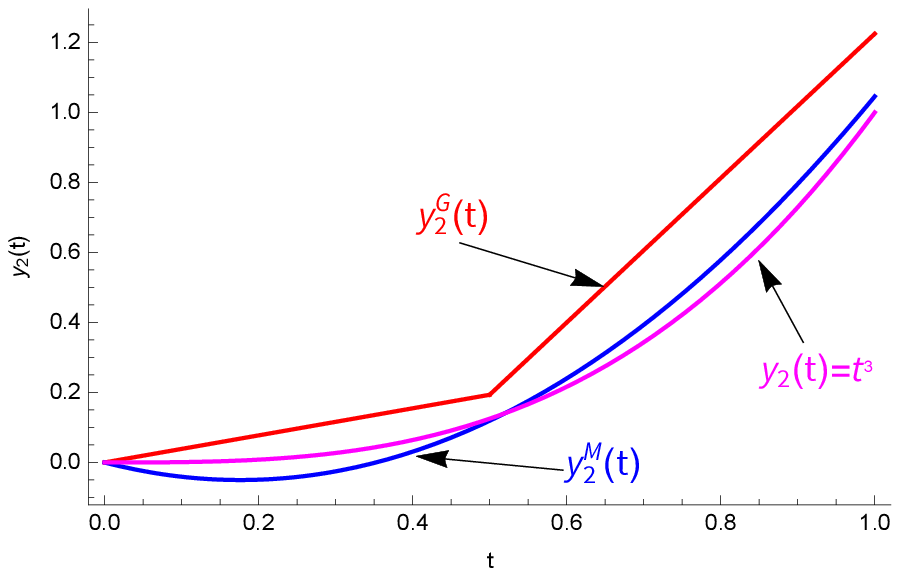}
\caption{Numerical results for problem \eqref{3.12} 
with $n=2$: results for the function $y_1(t)$ (left); 
results for the function $y_2(t)$ (right).}\label{fig:1}
\end{figure*}
\begin{table}[t]
\scriptsize
\centering
\caption{Numerical results for problem \eqref{3.12} 
with different values of $n$.}\label{tab:1}
\begin{tabular}{llllllllllllllll}
\hline
&\multicolumn{7}{c}{GHFs} && \multicolumn{7}{c}{MHFs}\\
\cline{2-8}\cline{10-16}
&\multicolumn{2}{c}{$y_1(t)$} &&\multicolumn{2}{c}{$y_2(t)$}
&&&&\multicolumn{2}{c}{$y_1(t)$} &&\multicolumn{2}{c}{$y_2(t)$}&\\
\cline{2-3}\cline{5-6}\cline{10-11}\cline{13-14}
$n$ & $e_{1,n}^G$ & $\rho_{1,n}^G$ & & $e_{2,n}^G$ & $\rho_{2,n}^G$ 
&  & CPU Time &&$e_{1,n}^M$ & $\rho_{1,n}^M$ & & $e_{2,n}^M$ 
& $\rho_{2,n}^M$ &  & CPU Time \\ \hline
$2$ &$1.68e-1$ & $2.32$&& $2.24e-1$ & $1.91$  && $0.000$  
&& $2.08e-2$&$4.48$&&$4.57e-2$&$3.99$&&$0.000$\\
$4$ &$3.37e-2$ & $2.13$&&$5.97e-2$  &$2.00$   
&&$0.000$   && $9.33e-4$&$3.45$&&$2.87e-3$&$3.69$&&$0.016$\\
$8$ &$7.68e-3$ & $2.02$&& $1.49e-2$  &$1.99$  &&$0.016$   
&& $8.51e-5$&$3.47$&&$2.23e-4$&$2.86$&&$0.016$\\
$16$ &$1.89e-3$& $1.98$&& $3.74e-3$  &$1.99$  &&$0.031$   
&& $7.69e-6$&$3.48$&&$3.08e-5$&$2.93$&&$0.047$\\
$32$ &$4.79e-4$ &$1.98$ && $9.43e-4$ &$1.99$  &&$0.172$   
&& $6.88e-7$&$3.49$&&$4.03e-6$&$2.96$&&$0.156$\\
$64$ &$1.21e-4$ & $1.98$&& $2.37e-4$  &$1.99$ &&$0.516$   
&& $6.12e-8$&$3.50$&&$5.16e-7$&$2.98$&&$0.859$\\
$128$ &$3.06e-5$ & $1.99$&& $5.97e-5$ &$2.00$ &&$4.453$   
&& $5.42e-9$&$3.50$&&$6.53e-8$&$2.99$&&$4.609$\\
$256$ &$7.70e-6$ & $2.00$&& $1.49e-5$ &$1.99$ &&$23.812$  
&& $4.80e-10$&$3.50$&&$8.21e-9$&$3.01$&&$23.641$\\
$512$ &$1.93e-6$ & ---&& $3.75e-6$ &---       
&&$183.688$ && $4.25e-11$ & ---&& $1.02e-9$ 
&---       &&$196.406$\\ \hline
\end{tabular}
\end{table}
\begin{figure*}
\centering
\includegraphics[scale=0.70]{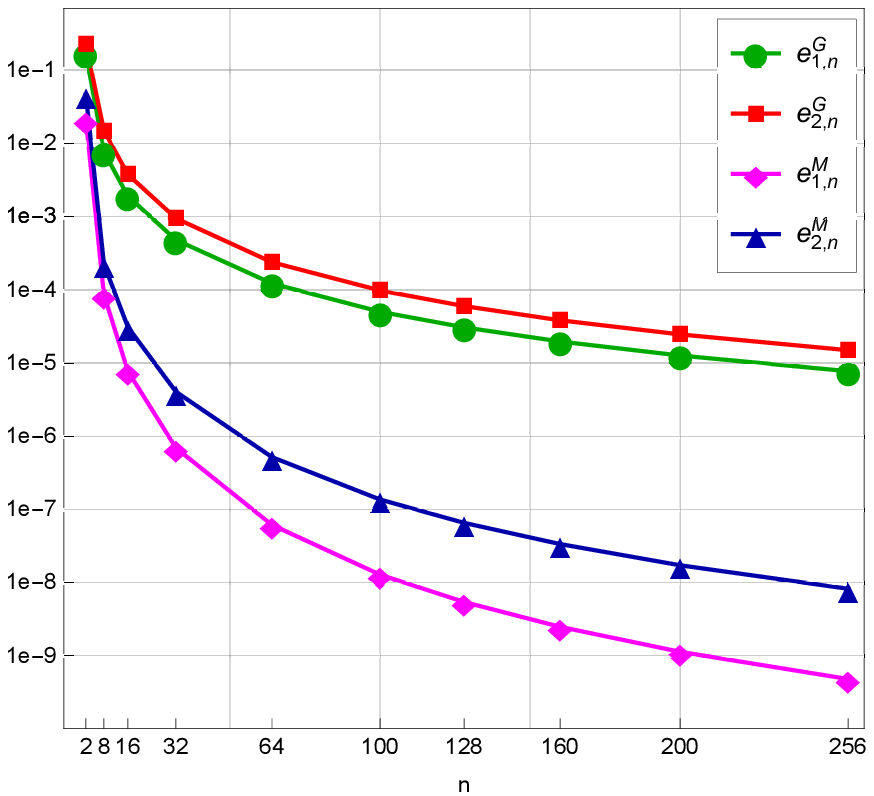}
\includegraphics[scale=0.70]{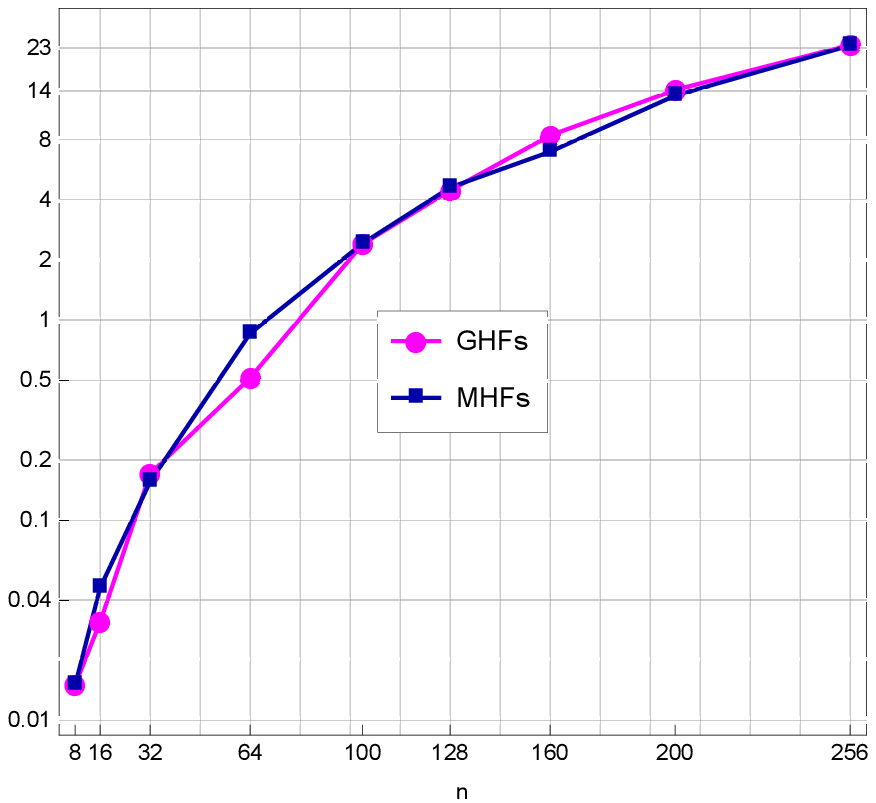}
\caption{Errors of the numerical method 
for solving problem \eqref{3.12} (left); 
CPU time in seconds (right).}\label{fig:2}
\end{figure*}
\end{example}

\vspace*{4cm}

\begin{example}
\label{ex2}
Consider the following system of linear FDEs:
\begin{equation}
\label{3.12-1}
\begin{cases}
{_0^CD_t^{\alpha}}y_1(t)=y_1(t)-2y_2(t)+4\cos(t)-2\sin(t),\quad   y_1(0)=1,\\
{_0^CD_t^{\alpha}}y_2(t)=3y_1(t)-4y_2(t)+5\cos(t)-5\sin(t),\quad  y_2(0)=2,
\end{cases}
\end{equation}
on the interval $0\leq t \leq 10$. The exact solution of this 
problem, when $\alpha=1$, is \cite{Atkinson}
\begin{equation*}
y_1(t)=\cos(t)+\sin(t),\quad y_2(t)=2 \cos(t).
\end{equation*}
We set $\alpha=1$ and solve the problem. By considering the same notations 
as introduced in Example~\ref{ex1}, we report the numerical results in 
Table~\ref{tab:2.1}. Furthermore, in Figures~\ref{fig:2.1} and \ref{fig:2.2}, 
the numerical solutions based on GHFs and MHFs, obtained by different values 
of $\alpha$ and $n=32$, together with the exact solution with $\alpha=1$, 
are plotted. As it could be expected, the numerical solution is close to the 
exact solution of the corresponding first order problem when $\alpha$ is close to $1$.
\begin{figure*}
\centering
\includegraphics[scale=0.70]{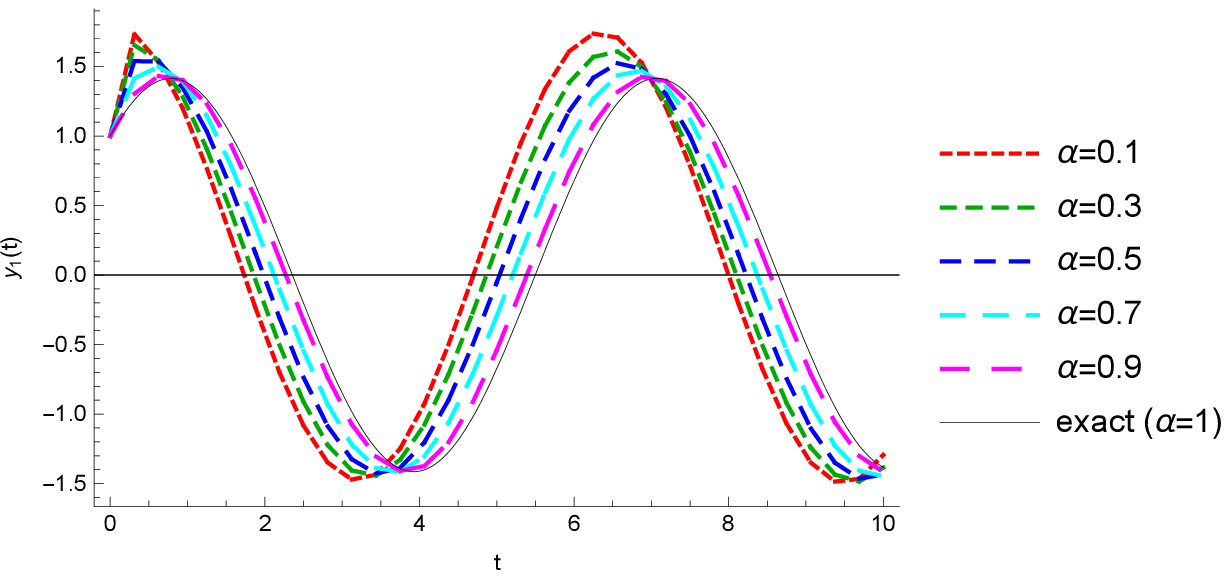}
\includegraphics[scale=0.70]{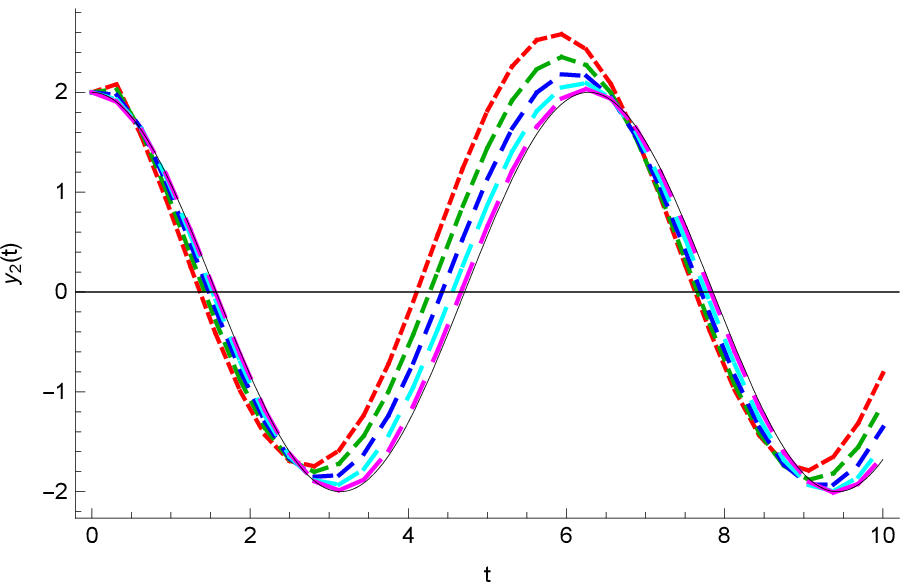}
\caption{Numerical results for problem \eqref{3.12-1} based on GHFs 
with different values of $\alpha$ and $n=32$: results for the function 
$y_1(t)$ (left); results for the function $y_2(t)$ (right).}\label{fig:2.1}
\end{figure*}
\begin{figure*}
\centering
\includegraphics[scale=0.70]{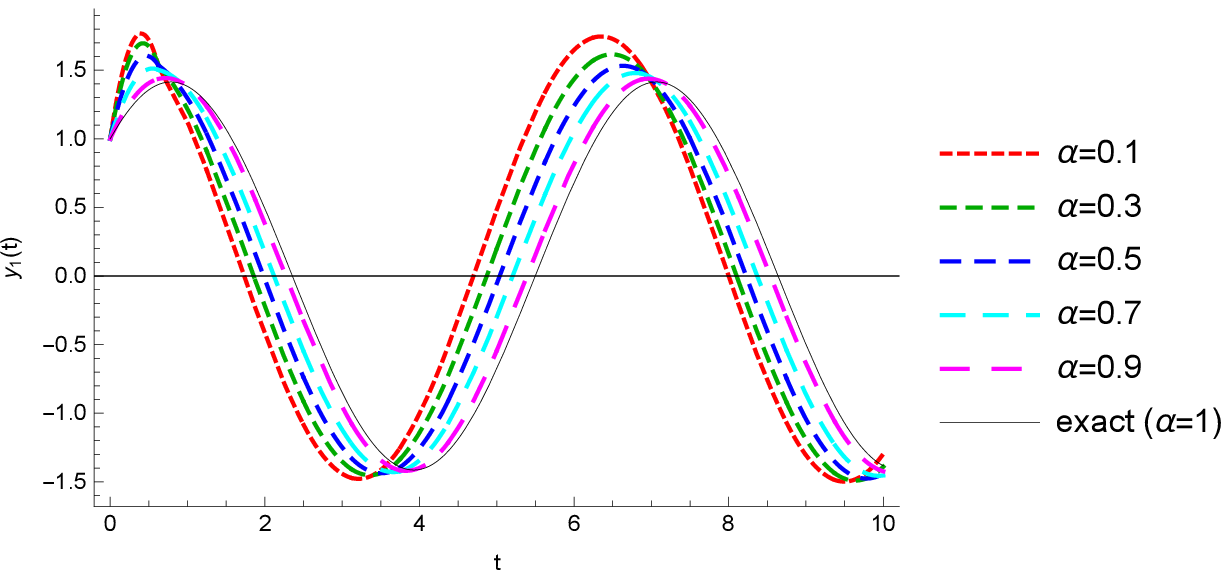}
\includegraphics[scale=0.70]{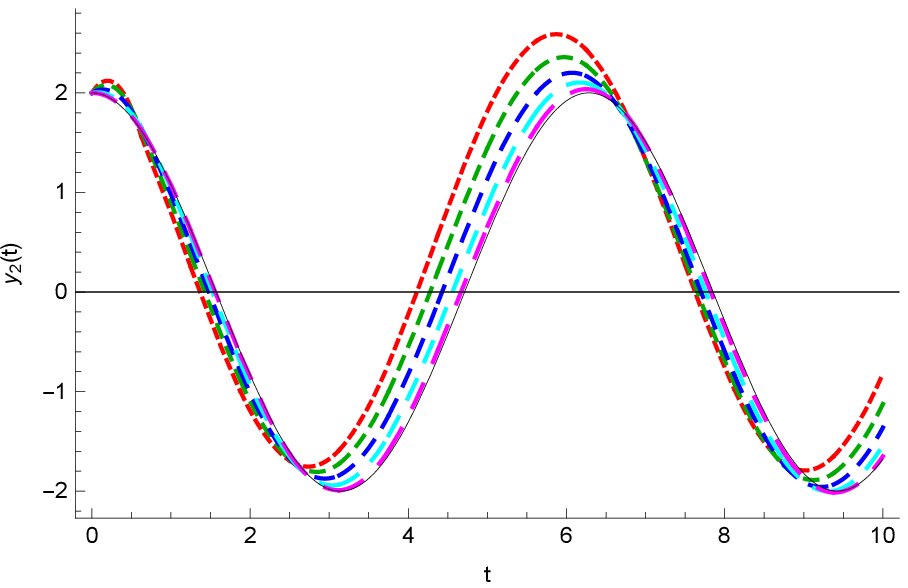}
\caption{Numerical results for problem \eqref{3.12-1} based on MHFs 
with different values of $\alpha$ and $n=32$: results for the function 
$y_1(t)$ (left); results for the function $y_2(t)$ (right).}\label{fig:2.2}
\end{figure*}
\end{example}


\section{Application to the SEIRS-$\alpha$ epidemic model}
\label{sec:4}

We now apply the numerical method introduced 
in Section~\ref{sec:3} to a nonlocal fractional 
order SEIRS mathematical model, which was 
recently proposed in \cite{Rosa}.


\subsection{Description of the model}

Human respiratory syncytial virus (HRSV) is a virus that causes respiratory 
tract infections. We refer a reader interested in this virus to 
\cite{Rosa} and references therein. Here we just mention that HRSV
is a principal cause of lower respiratory tract infections and hospital visits 
during infancy and childhood. There is an annual epidemic in temperate climates 
during the winter season, while in tropical climates this infection is most common 
throughout the rainy season. Since protective immunity is induced by natural 
infection with HRSV more than many other respiratory viral infections, people 
can be infected multiple times even within a single HRSV season. 

A mathematical model can show how the infectious disease with HRSV 
progresses and what are the outcomes of an epidemic of this virus. Recently, 
a compartmental model was proposed in \cite{Rosa}, 
based on stratifying the population into four health states: susceptible to the 
infection, denoted by $S$; a group of individuals $E$ who have been infected but are 
not infectious yet, which become infectious at a rate $\varepsilon$;
infected and infectious, denoted by $I$; and recovered individuals $R$. 
A particular property of HRSV is that immunity after infection is temporary, i.e., 
the recovered individuals become susceptible again \cite{Weber}, hence, 
the model is called a SEIRS model. The authors of \cite{Rosa} considered that 
the annual recruitment rate is seasonal due to schools opening/closing 
periods and proposed the following system of FDEs:
\begin{equation}
\label{4.1}
\begin{cases}
{_0^CD_t^{\alpha}} S(t)=\lambda(t)-\mu S(t)-\beta(t)S(t)I(t)+\gamma R(t),\\
{_0^CD_t^{\alpha}} E(t)=\beta(t)S(t)I(t)-\mu E(t)-\varepsilon E(t),\\
{_0^CD_t^{\alpha}} I(t)=\varepsilon E(t)-\mu I(t)-\nu I(t)\\
{_0^CD_t^{\alpha}} R(t)=\nu I(t)-\mu R(t)-\gamma R(t),
\end{cases}
\end{equation}
with given initial conditions
\begin{equation*}
S(0),E(0),I(0),R(0)\geq 0,
\end{equation*}
where $\mu$ denotes the birth rate, which was assumed equal to the mortality 
rate, $\gamma$ is the rate of loss of immunity, $\nu$ is the rate of loss 
of infectiousness, $\beta$ denotes the transmission parameter, which is modeled 
by the cosine function as $\beta(t)=b_0(1+b_1\cos(2\pi t+\Phi))$, in which $b_0$ 
is the mean of $\beta$ and $b_1$ is the amplitude of the seasonal fluctuation, 
$\lambda(t)=\mu(1+c_1\cos(2\pi t+\Phi))$ is the recruitment rate, which includes 
newborns and immigrants, with $c_1$ as the amplitude of the seasonal fluctuation, 
and where $_0^CD_t^{\alpha}$ denotes the left Caputo derivative of order $\alpha\in (0,1]$. 
In the parameters $\beta$ and $\lambda$, the parameter $\Phi$ is an angle that is 
chosen in agreement with real data. Note that by introducing the group $E$, a 
latency period is included in the model, which is assumed equal to the time 
between infection and the first symptoms.

\begin{table}[t]
\scriptsize
\centering
\caption{Numerical results for problem \eqref{3.12-1} 
with different values of $n$ and $\alpha=1$.}\label{tab:2.1}
\begin{tabular}{llllllllllllllll}
\hline
&\multicolumn{7}{c}{GHFs} && \multicolumn{7}{c}{MHFs}\\
\cline{2-8}\cline{10-16}
&\multicolumn{2}{c}{$y_1(t)$} &&\multicolumn{2}{c}{$y_2(t)$}
&&&&\multicolumn{2}{c}{$y_1(t)$} &&\multicolumn{2}{c}{$y_2(t)$}&\\
\cline{2-3}\cline{5-6}\cline{10-11}\cline{13-14}
$n$ & $e_{1,n}^G$ & $\rho_{1,n}^G$ & & $e_{2,n}^G$ & $\rho_{2,n}^G$ 
&  & CPU Time &&$e_{1,n}^M$ & $\rho_{1,n}^M$ & & $e_{2,n}^M$ 
& $\rho_{2,n}^M$ &  & CPU Time \\ \hline
$2$ &$2.33e+0$ & $2.22$&& $2.15e+0$ & $3.01$  && $0.000$  
&& $2.03e+0$&$1.26$&&$1.87e+0$&$1.53$&&$0.000$\\
$4$ &$5.01e-1$ & $1.93$&&$2.67e-1$  &$2.31$   
&&$0.031$   && $8.47e-1$&$3.23$&&$6.49e-1$&$2.48$&&$0.000$\\
$8$ &$1.31e-1$ & $2.04$&& $5.38e-2$  &$1.86$  &&$0.031$   
&& $9.00e-2$&$3.67$&&$1.16e-1$&$3.94$&&$0.016$\\
$16$ &$3.18e-2$& $1.94$&& $1.48e-2$  &$1.95$  &&$0.031$   
&& $7.05e-3$&$4.07$&&$7.56e-3$&$3.87$&&$0.016$\\
$32$ &$8.29e-3$ &$2.00$ && $3.82e-3$ &$1.99$  &&$0.031$   
&& $4.20e-4$&$4.04$&&$9.63e-4$&$3.72$&&$0.031$\\
$64$ &$2.07e-3$ & $2.00$&& $9.63e-4$  &$2.00$ &&$0.078$   
&& $2.55e-5$&$4.00$&&$3.91e-5$&$3.83$&&$0.109$\\
$128$ &$5.17e-4$ & $2.00$&& $2.41e-4$ &$2.00$ &&$0.344$   
&& $1.59e-6$&$4.00$&&$2.74e-6$&$3.91$&&$0.344$\\
$256$ &$1.29e-4$ & $2.00$&& $6.03e-5$ &$2.01$ &&$1.094$  
&& $9.93e-8$&$4.00$&&$1.82e-7$&$3.96$&&$1.172$\\
$512$ &$3.23e-5$ & ---&& $1.50e-5$ &---       
&&$3.266$ && $6.20e-9$ & ---&& $1.17e-8$ 
&---       &&$3.656$\\ \hline
\end{tabular}
\end{table}


\vspace*{4cm}

\subsection{Numerical results}

\begin{figure*}
\centering
\includegraphics[scale=0.70]{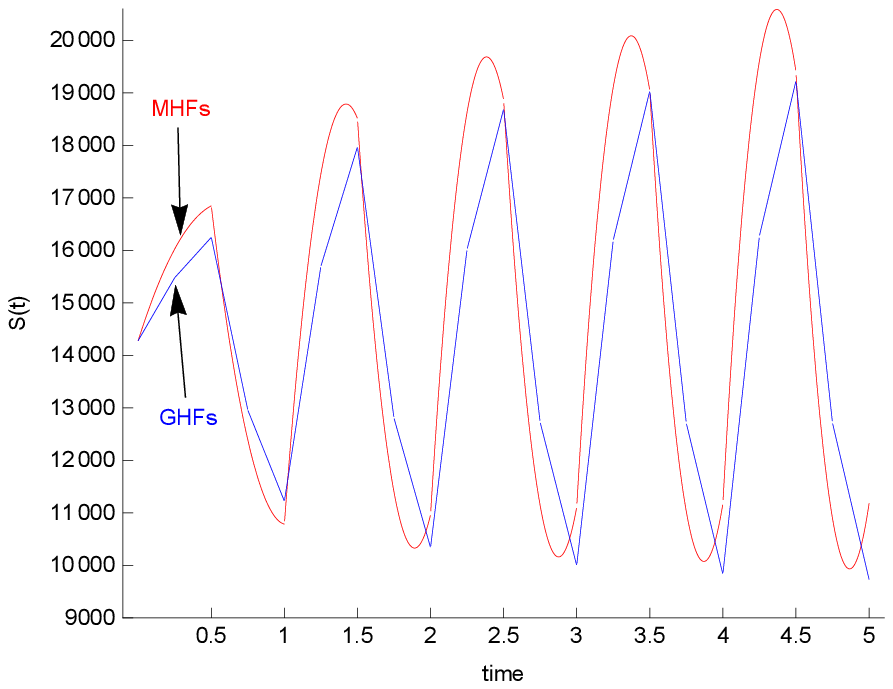}
\includegraphics[scale=0.70]{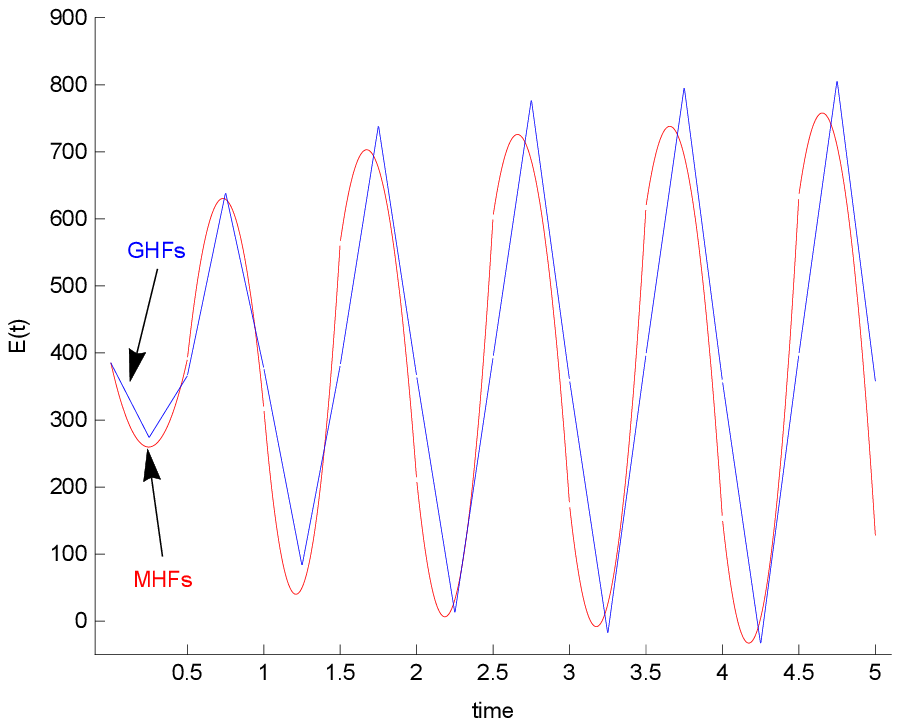}
\includegraphics[scale=0.70]{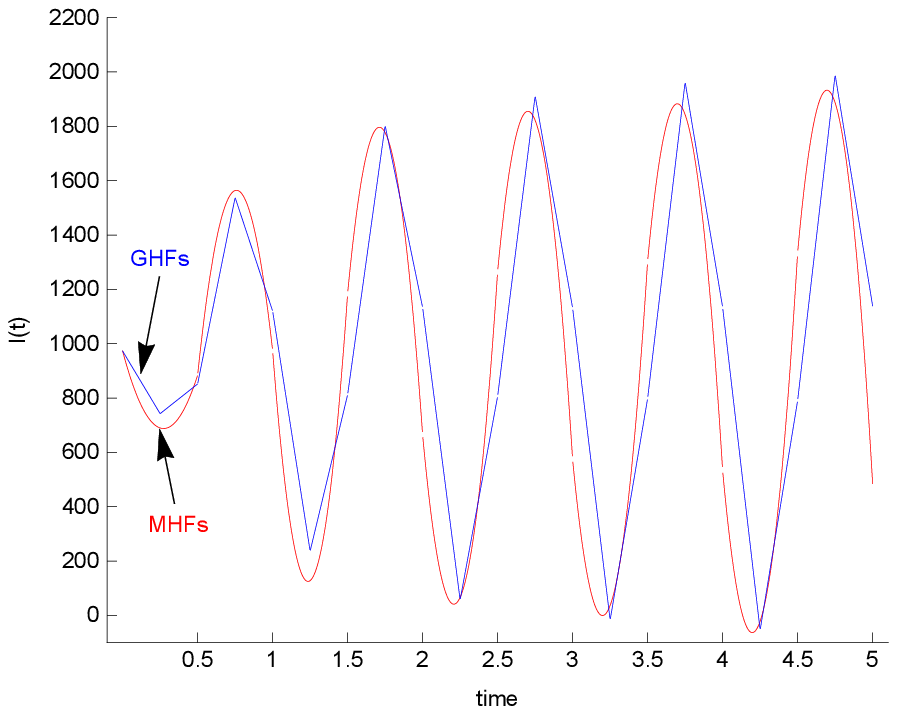}
\includegraphics[scale=0.70]{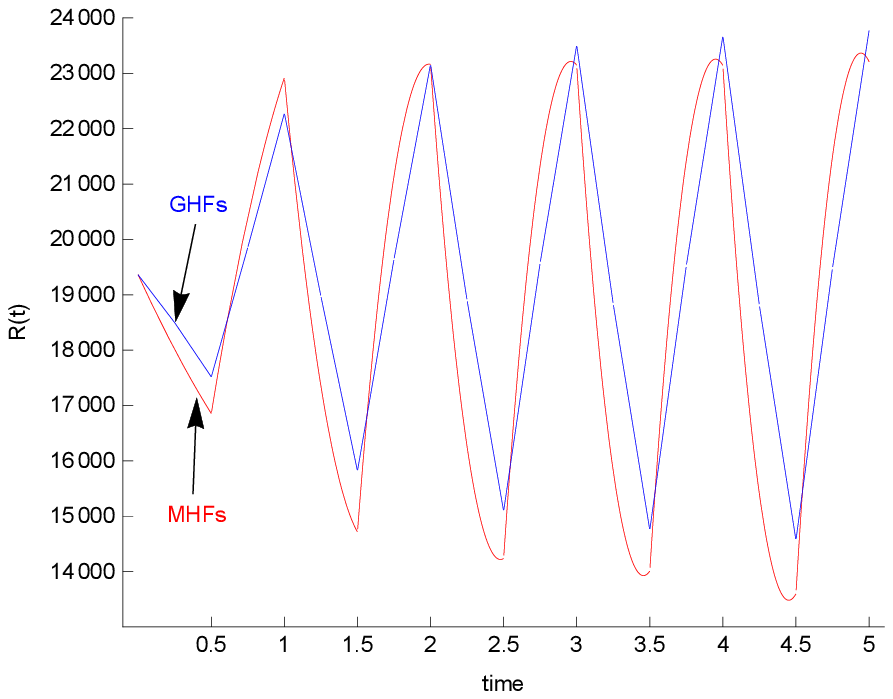}
\caption{Numerical results of the state variables of the SEIRS-$\alpha$ 
model given by \eqref{4.1} based on GHFs and MHFs with $n=20$, 
considering $\alpha=0.993$.}\label{fig:3}
\end{figure*}

\begin{figure*}
\centering
\includegraphics[scale=0.70]{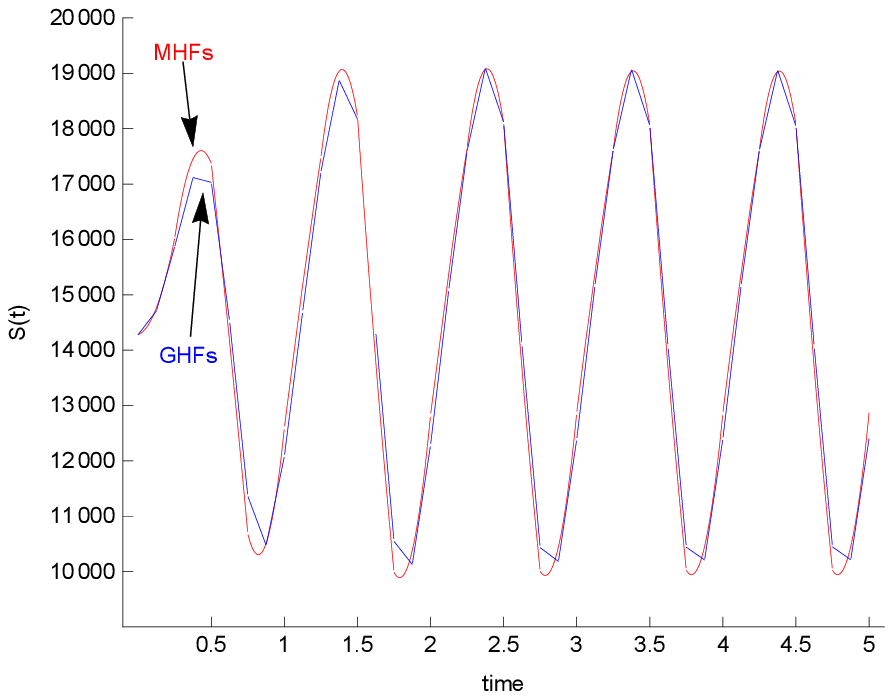}
\includegraphics[scale=0.70]{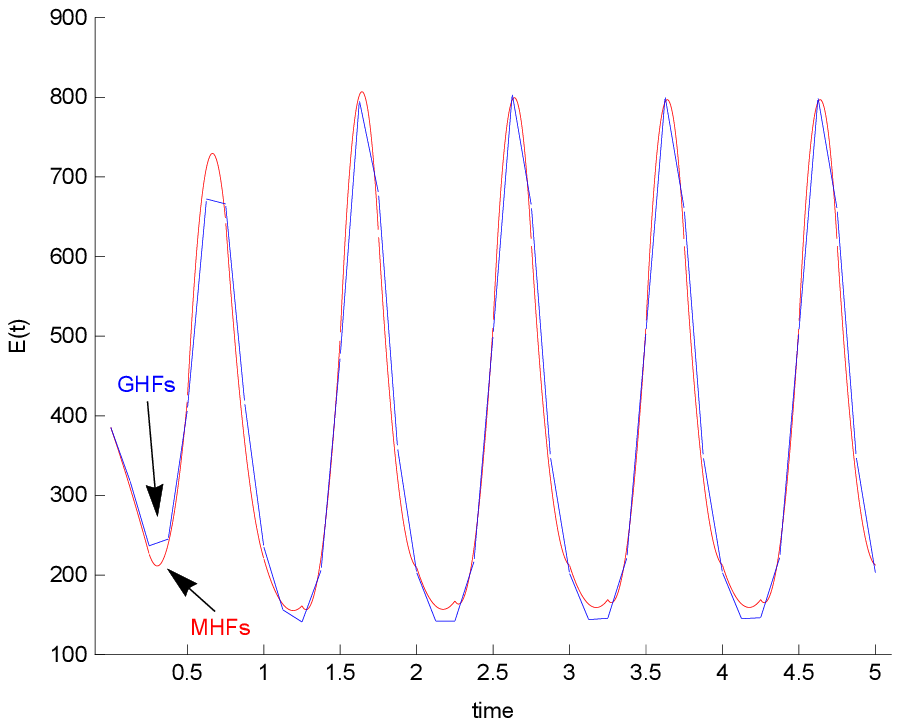}
\includegraphics[scale=0.70]{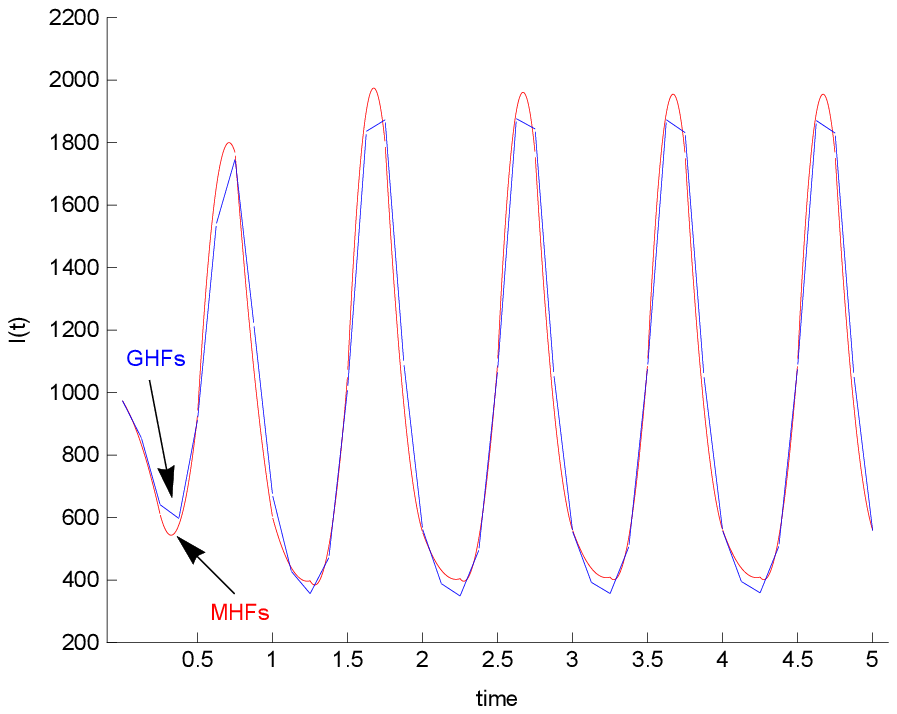}
\includegraphics[scale=0.70]{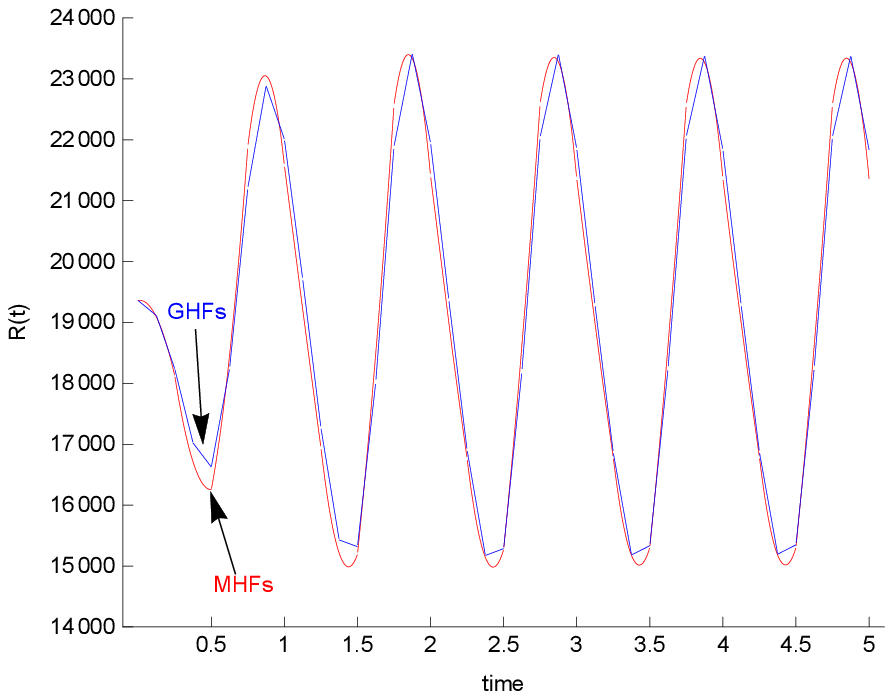}
\caption{Numerical results of the state variables of the SEIRS-$\alpha$ 
model given by \eqref{4.1} based on GHFs and MHFs with $n=40$, 
considering $\alpha=0.993$.}\label{fig:4}
\end{figure*}

\begin{figure*}
\centering
\includegraphics[scale=0.70]{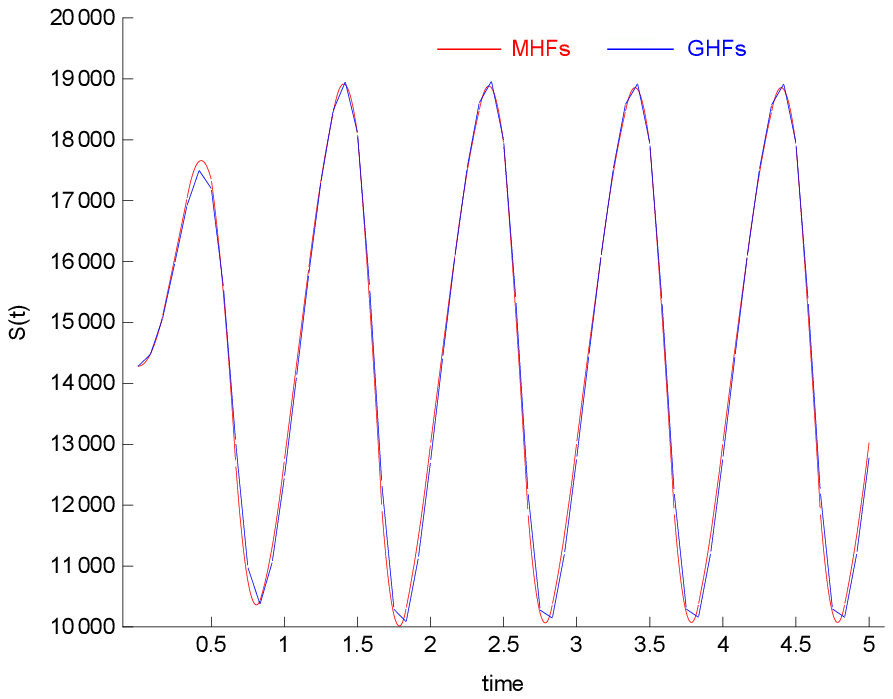}
\includegraphics[scale=0.70]{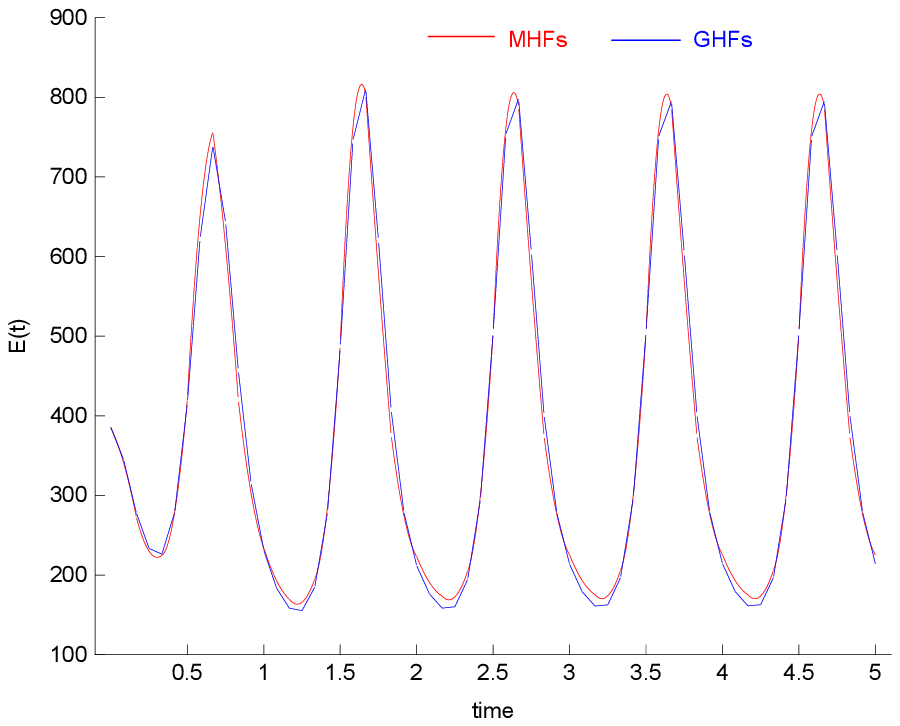}
\includegraphics[scale=0.70]{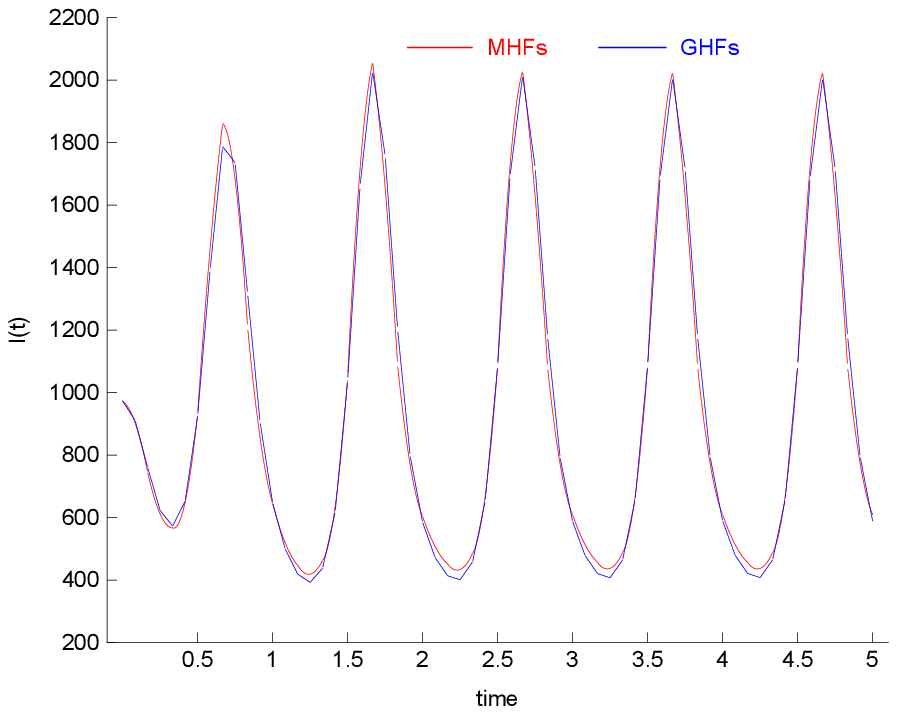}
\includegraphics[scale=0.70]{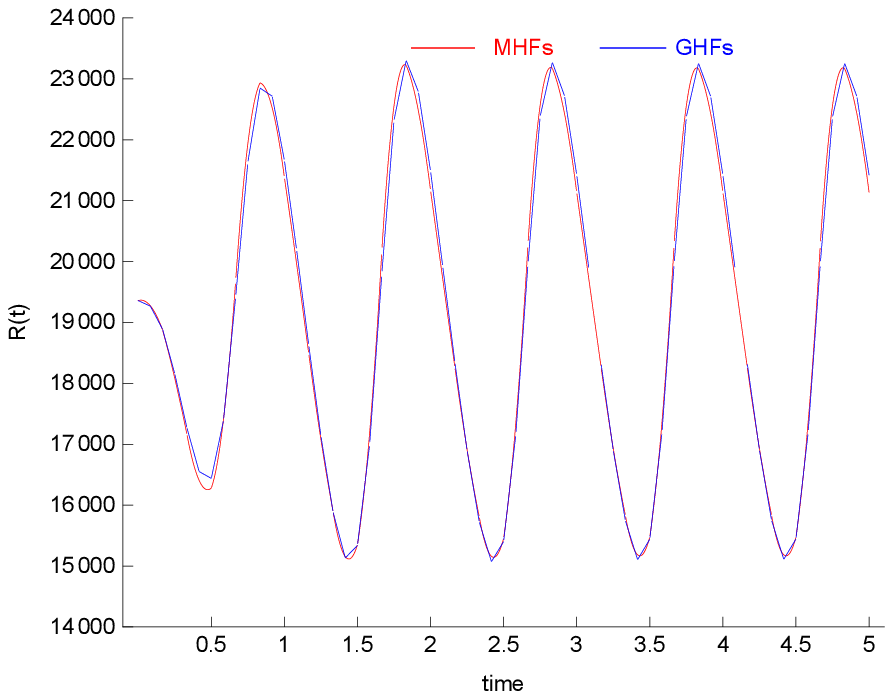}
\caption{Numerical results of the state variables of the SEIRS-$\alpha$ 
model given by \eqref{4.1} based on GHFs and MHFs with $n=60$, 
considering $\alpha=0.993$.}\label{fig:5}
\end{figure*}

\begin{figure*}
\centering
\includegraphics[scale=0.70]{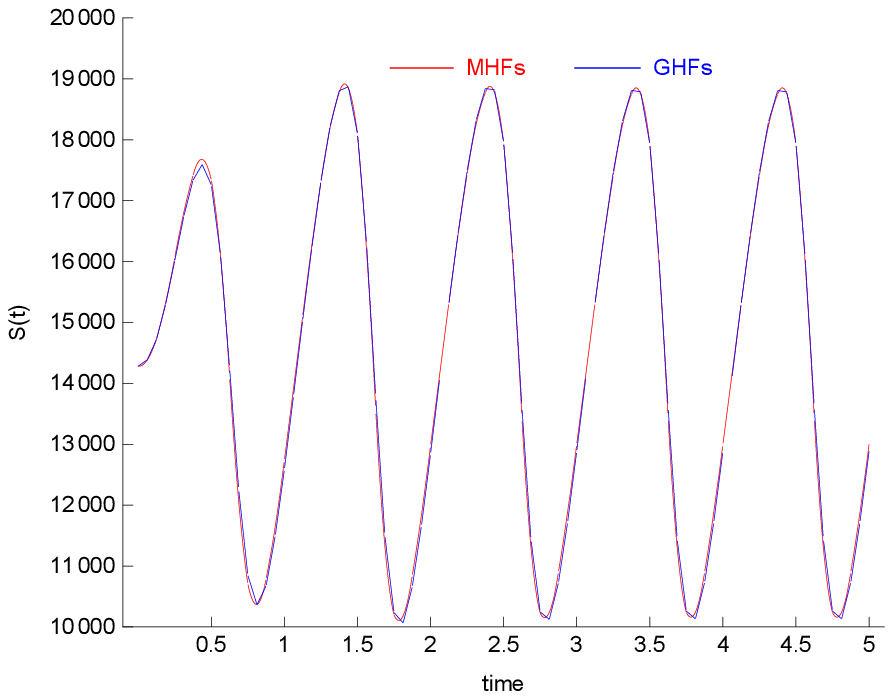}
\includegraphics[scale=0.70]{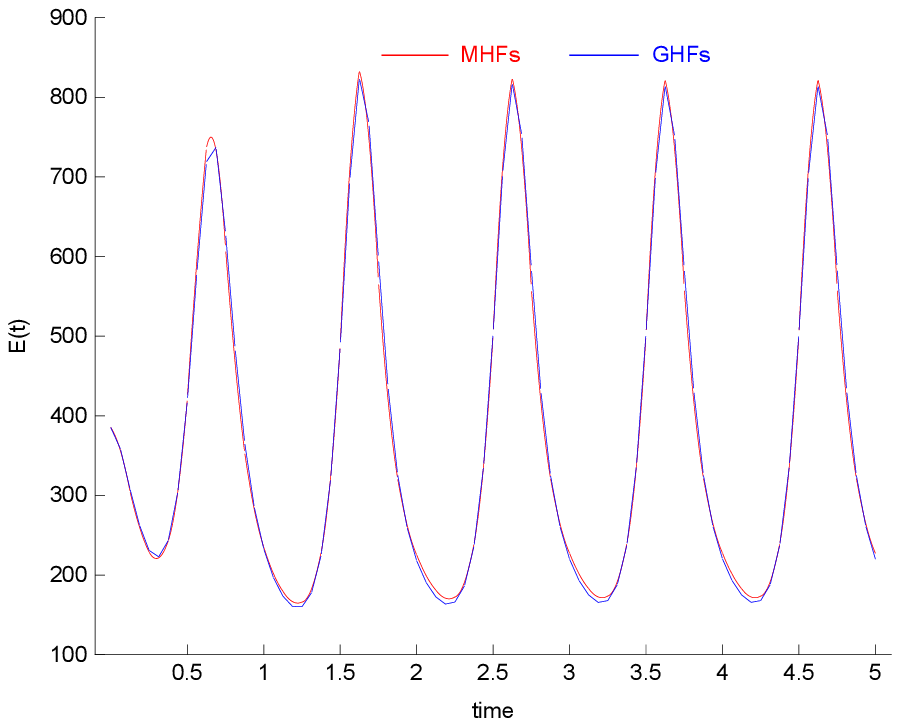}
\includegraphics[scale=0.70]{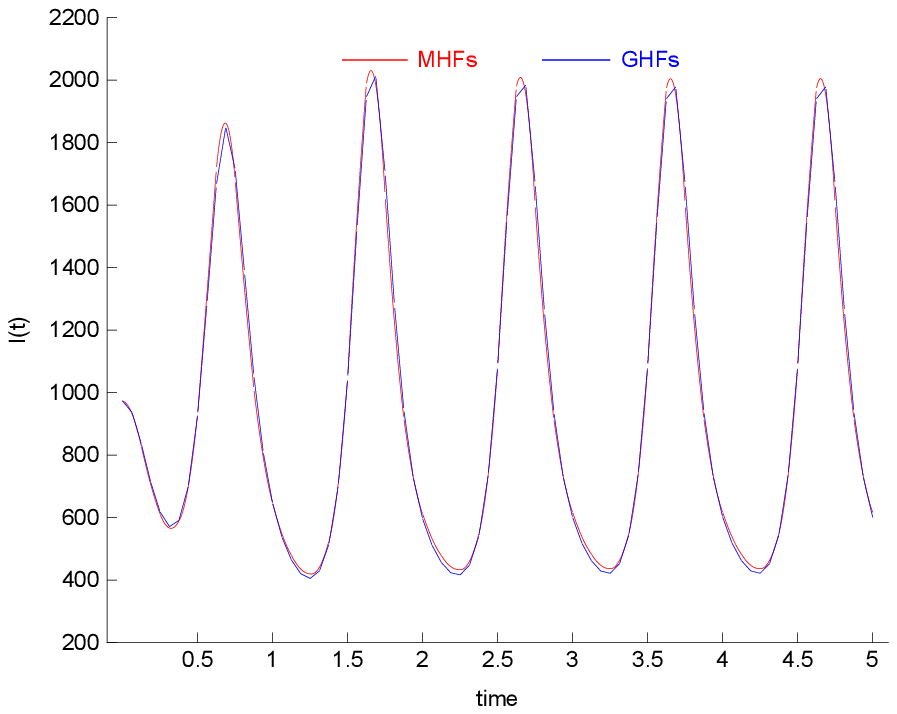}
\includegraphics[scale=0.70]{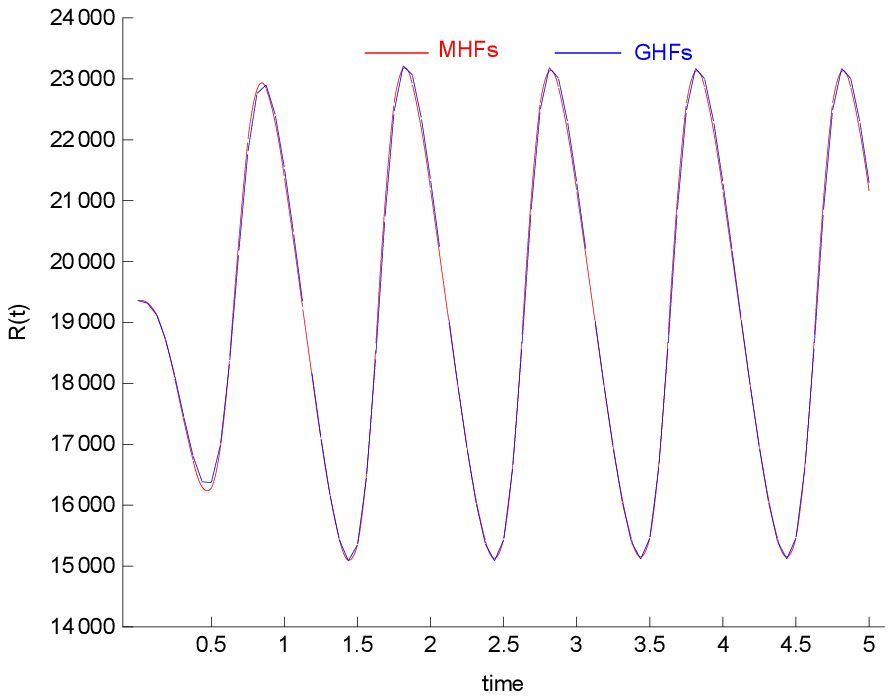}
\caption{Numerical results of the state variables of the SEIRS-$\alpha$ 
model given by \eqref{4.1} based on GHFs and MHFs with $n=80$, 
considering $\alpha=0.993$.}\label{fig:6}
\end{figure*}

Based on data obtained from the Florida Department of Health, authors 
in \cite{Rosa} searched the fractional order of differentiation, $\alpha$, 
that best fits the data on the reported number of positive tests of HRSV disease, 
per month, during 35 months, precisely between September 2011 and July 2014 in the state 
of Florida (excluding North region). They found that with $\alpha=0.993$ the model 
fits quite well the data of HRSV disease. Using the values of the parameters $\mu$, 
$\nu$, $\gamma$, $\varepsilon$, $b_0$, $b_1$, $c_1$, and $\Phi$ as given in 
Table~\ref{tab:2}, borrowed from \cite{Rosa}, and the initial conditions given 
in Table~\ref{tab:3}, we applied our numerical scheme to solve the SEIRS-$\alpha$ 
model given by \eqref{4.1} with $\alpha=0.993$ and different values of $n$. 
Since World Health Organization goals are usually fixed for five years periods 
for most diseases, it was assumed that $\tau=5$. By considering $n=20, 40, 60, 80$, 
the numerical results of the state variables $S(t)$, $E(t)$, $I(t)$ and $R(t)$, 
obtained by the proposed method based on GHFs and MHFs, are reported 
in Figures~\ref{fig:3}--\ref{fig:6}. From these figures, we see that by 
increasing the value of $n$, the results are in full agreement 
with those of \cite{Rosa}, which in contrast with our approach
are obtained by indirect methods.
\begin{table}[h]
\scriptsize
\centering
\caption{Parameters of the SEIRS-$\alpha$ model borrowed from \cite{Rosa}.}
\label{tab:2}
\begin{tabular}{llllllll}
\hline
$\mu$   & $\nu$ & $\gamma$  & $\varepsilon$  & $b_0$   &  $b_1$   &  $c_1$ & $\Phi$\\
\hline
$0.0113$ & $36$ & $1.8$ & $91$ & $88.25$ & $0.17$  & $0.17$ & $\frac{\pi}{2}$\\
\hline
\end{tabular}
\end{table}
\begin{table}[h]
\scriptsize
\centering
\caption{Initial conditions in terms of percentage of total population.}
\label{tab:3}
\begin{tabular}{llll}
\hline
$S(0)$   & $E(0)$ &  $I(0)$  & $R(0)$ \\
\hline
$0.4081$ & $0.0110$ & $0.0278$ & $0.5531$ \\
\hline
\end{tabular}
\end{table}


\section{Concluding remarks}
\label{sec:5}

We introduced a new numerical approach for solving systems 
of fractional differential equations (FDEs). In our method, 
two classes of hat functions, namely generalized hat functions (GHFs) 
and modified hat functions (MHFs), have been used as the basis functions. 
In this new scheme, the considered system of FDEs is easily reduced 
to a system of nonlinear algebraic equations that, according to the 
computational complexity, needs few computational efforts 
in both cases of basis functions. By applying the method to two test problems, 
we see that the numerical results confirm the $O(h^2)$ accuracy order of the GHFs 
and the $O(h^3)$ accuracy order of the the MHFs. Finally, the method is successfully 
applied to a mathematical model in epidemiology given by a system of FDEs 
with application to human respiratory syncytial virus infection.


\begin{acknowledgements}
Torres was supported by the Portuguese national funding 
agency for science, research and technology (FCT),
within project UID/MAT/04106/2019 (CIDMA).
The authors are grateful to two anonymous referees for 
several positive and constructive comments, 
which helped them to improve the manuscript.
\end{acknowledgements}


\section*{Compliance with ethical standards}

\textbf{Conflict of interest}.
The authors declare that they have no conflict of interest.

\medskip

\noindent \textbf{Ethical approval}.
This article does not contain any studies with human
participants or animals performed by any of the authors.



\end{document}